\documentclass[oneside,english]{amsart}
\usepackage[T1]{fontenc}
\usepackage[latin9]{inputenc}
\usepackage{geometry}
\usepackage{float}
\usepackage{amsthm}
\usepackage{amstext}
\usepackage{amssymb}
\usepackage{graphicx}
\usepackage{esint}

\makeatletter

\numberwithin{equation}{section}
\numberwithin{figure}{section}
\theoremstyle{plain}
\newtheorem{thm}{\protect\theoremname}
  \theoremstyle{plain}
  \newtheorem{lem}[thm]{\protect\lemmaname}
  \theoremstyle{plain}
  \newtheorem{prop}[thm]{\protect\propositionname}
  \theoremstyle{remark}
  \newtheorem{rem}[thm]{\protect\remarkname}

\usepackage{color}

\makeatother

\usepackage{babel}
  \providecommand{\lemmaname}{Lemma}
  \providecommand{\propositionname}{Proposition}
  \providecommand{\remarkname}{Remark}
\providecommand{\theoremname}{Theorem}

\begin{document}
\global\long\def\intr{\int_{R}}
 \global\long\def\cbr#1{\left\{  #1\right\}  }
\global\long\def\rbr#1{\left(#1\right)}
\global\long\def\ev#1{\mathbb{E}{#1}}
\global\long\def\R{\mathbb{R}}
\global\long\def\norm#1#2#3{\Vert#1\Vert_{#2}^{#3}}
\global\long\def\pr#1{\mathbb{P}\rbr{#1}}
\global\long\def\cleq{\lesssim}
\global\long\def\ceq{\eqsim}
\global\long\def\conv{\rightarrow}
\global\long\def\Var#1{\text{Var}(#1)}
\global\long\def\TDD#1{{\color{red}TDD(#1)}}
\global\long\def\dd#1{\textnormal{d}#1}
\global\long\def\inti{\int_{0}^{\infty}}
\global\long\def\eqdef{:=}
\global\long\def\crr{\mathcal{C}([0,\infty),\R)}
\global\long\def\sb#1{\langle#1\rangle}
\global\long\def\pm#1{d_{P}\rbr{#1}}
\global\long\def\crt{\mathcal{C}([0,T],\R)}
\global\long\def\ab{[a;b]}
\global\long\def\TTV#1#2#3{\text{TV}^{#3}\!\rbr{#1,#2}}
\global\long\def\UTV#1#2#3{\text{UTV}^{#3}\!\rbr{#1,#2}}
\global\long\def\DTV#1#2#3{\text{DTV}^{#3}\!\rbr{#1,#2}}

\title{Exact representation of truncated variation of Brownian motion}

\author{Piotr Mi\l{}o\'{s}}

\address{Faculty of Mathematics, Informatics and Mechanics, University of
Warsaw\\
 Banacha 2, 02-097 Warszawa, Poland}

\email{pmilos@mimuw.edu.pl}

\date{11.11.2013}
\begin{abstract}
In the recent papers \cite{ochowski:2011fk,Lochowski:2013yq,ochowski:2013lr}
the truncated variation has been introduced, characterized and studied
in various stochastic settings. In this note we uncover an intimate
link to the Skorokhod problem. Further, we exploit it to give an explicit
representation of the truncated variation of a Brownian motion. More
precisely, we prove that the inverse of this process is, up to a minor
time shift, a Lévy subordinator with the exponent {\normalsize $\sqrt{2q}\tanh(c\sqrt{q/2})$. }{\normalsize \par}

This also gives a representation of a solution of the two-sided Skorokhod
problem for a Brownian motion. 
\end{abstract}
\maketitle

\section{Representation of the truncated variation}

Given $f:\ab\mapsto\R$, a càdlàg function, its truncated variation
on the interval $\ab$ is defined by 
\begin{equation}
\TTV f{\ab}c:=\sup_{n}\sup_{a\leq t_{1}<t_{2}<\ldots<t_{n}\leq b}\sum_{i=1}^{n-1}\phi_{c}\left(f(t_{i+1})-f(t_{i})\right),\quad c\geq0,\label{eq:tvDef}
\end{equation}
$\phi_{c}\left(x\right)=\max\left\{ |x|-c,0\right\} $ (note that
$\TTV f{\ab}0$ is the total variation). Closely related notions of
the upward truncated variation and the downward truncated variation
denoted by $\text{UTV}^{c}$, respectively $\text{DTV}^{c}$ are defined
by putting $\phi_{c}(x)=\max\left\{ x-c,0\right\} $, respectively
$\phi_{c}(x)=\max\left\{ -x-c,0\right\} $ in \eqref{eq:tvDef}. The
truncated variation has two interesting variational characterizations
- \cite[(1.2) and (2.2)]{Lochowski:2013yq}. We recall one of them:
\begin{equation}
\TTV f{\ab}c=\inf\cbr{\TTV g{\ab}0:g\,\text{such that }\norm{g-f}{\infty}{}\leq\frac{1}{2}c},\label{eq:variational}
\end{equation}
where $\norm g{\infty}{}:=\sup\cbr{|g(x)|:x\in\ab}$. Moreover, the
infimum is attained for some function $g^{c}$, which can be effectively
characterized. We skip the precise description at the moment (instead
we refer the reader to \cite[Section 2.1]{Lochowski:2013yq} and the
proof of Proposition \ref{prop:skorohodEmbending} below), yet the
basic idea is simple. Following \cite[Remark 2.4]{ochowski:2011lr}
we notice that {}``$g^{c}$ is the most lazy function, which changes
its value only if it is necessary to stay in the tube defined by $\norm{g^{c}-f}{\infty}{}\leq c/2$''.
This can be seen on the following picture: 
\begin{figure}[H]
\caption{An example of function (in red) and $g^{c}$ (in various colors).}

\includegraphics[scale=0.45]{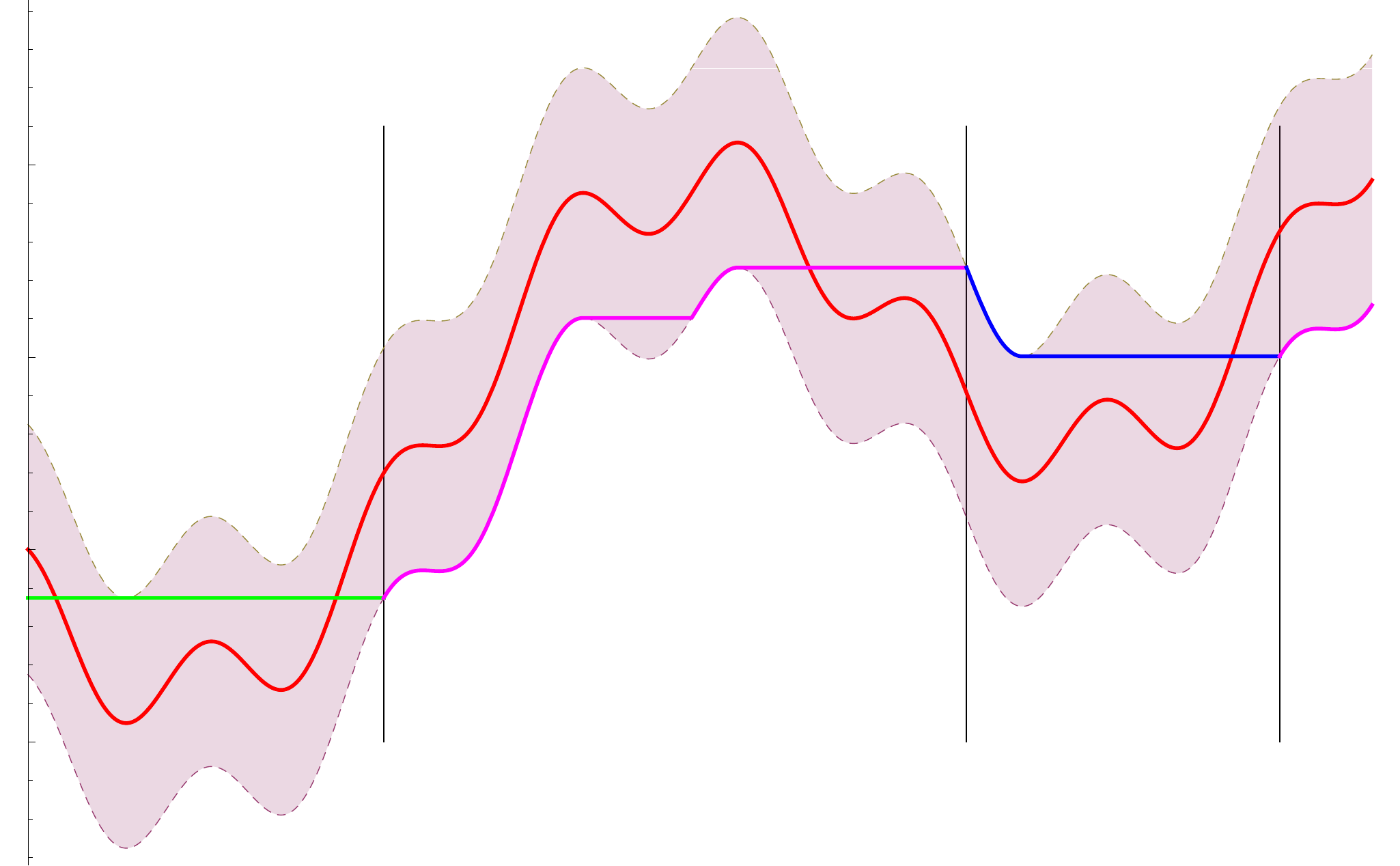}
\end{figure}

This observation, explored from a different perspective, will lead
us to a new characterization of $g^{c}$ in the language of the Skorokhod
problem. First we recall that 
\begin{lem}
\label{lem:skorohod}Let $f\in\mathcal{C}(\ab,\R)$, $c>0$. For any
$x\in[f(a)-c/2,f(a)+c/2]$ there exist unique functions $f_{1},f_{2}\in\mathcal{C}(\ab,\R)$
such that
\begin{enumerate}
\item $f_{1}(a)=f_{2}(a)=0$,
\item $h(s)\eqdef f(s)-\rbr{x+f_{1}(s)-f_{2}(s)}\in[-c/2,c/2]$ for any
$s\in\ab$,
\item $f_{1}$ is increasing, continuous, and the corresponding measure
$\dd{f_{1}}$ is carried by $\cbr{t:h(t)=-c/2}$,
\item $f_{2}$ is increasing, continuous, and the corresponding measure
$\dd{f_{2}}$ is carried by $\cbr{t:h(t)=c/2}$.
\end{enumerate}
\end{lem}
The lemma is not a new result. Its roots might be traced back to the
seminal works of Skorokhod as defined in \cite[Lemma VI.2.1]{Revuz:1991kx}.
A more general versions of the problem were studied later on, culminating
in \cite{Burdzy:2009lr} (see \cite[Theorem 2.6]{Burdzy:2009lr} for
the uniqueness and existence result). Now we are ready to link the
truncated variation and the Skorokhod problem. To this end we define
$\tau_{U}\eqdef\inf\cbr{s\geq a:\rbr{f(s)-\inf_{a\leq u\leq s}f(u)}\geq c}$
and $\tau_{D}\eqdef\inf\cbr{s\geq a:\rbr{\sup_{a\leq u\leq s}f(u)-f(s)}\geq c}$.
We put
\[
x\eqdef\begin{cases}
f(\tau_{U})-c/2, & \text{ if }\tau_{U}<\tau_{D},\\
f(\tau_{D})+c/2, & \text{ if }\tau_{U}>\tau_{D}.
\end{cases}
\]

\begin{prop}
\label{prop:skorohodEmbending}For this special choice of $x$ the
solution of the Skorokhod problem given by Lemma \ref{lem:skorohod}
fulfills 
\[
g^{c}(t)=x+f_{1}(s)-f_{2}(s).
\]
Consequently for any $t\in[a;b]$ we have
\[
\TTV f{[a;t]}c=f_{1}(t)+f_{2}(t).
\]
Moreover, one easily observes that $f_{1}(s)=f_{2}(s)=0$ if $s\in[a;\tau_{U}\wedge\tau_{D}]$. \end{prop}
\begin{rem}
The choice of $x$ for the starting point might seem little mysterious.
Again this is a sign of the {}``laziness''. One chooses $x$ such
that $g^{c}$ can remain constant as long as it possible (which can
be seen on the picture above, in which case we have $\tau_{U}<\tau_{D}$).
\end{rem}

\begin{proof}
This proposition is an easy consequence of known results thus we present
it in a sketchy way. We are going to check that $f_{1}$ is the same
as $\text{UTV}^{c}$ and $f_{2}$ is $\text{DTV}^{c}$. By \cite[(2.4)]{Lochowski:2013yq},
we know that the unique minimizer, \textbf{$g^{c}$,} in \eqref{eq:variational}
is given by
\[
g^{c}(s)\eqdef x+\UTV f{[a;s]}c-\DTV f{[a;s]}c.
\]
By Lemma \ref{lem:skorohod} it is enough to check conditions 1-4.
Condition 1 is trivial, condition 2 holds by \cite[(2.1)]{Lochowski:2013yq}.
Concerning condition $3$, it is easy to check that $f_{1}$ is continuous
and increasing. We now recall \cite[Section 2.2]{ochowski:2011lr},
it studies function $f^{c}$, which by \cite[Theorem 2.1]{ochowski:2011lr}
is the same as our $g^{c}$. We see this function grows if and only
if the running maximum $M_{k}^{c}(s)$ grows and $s\in[T_{U,k}^{c};T_{D,k}^{c})$.
It is straightforward to check that these are precisely $s$ for which
$h(s)=-c/2$.

\end{proof}

\section{Truncated variation of Brownian Motion}

We will apply the results of the previous section to study the truncated
variation of a Brownian motion. Theorem \ref{thm:localTimeRepresentation}
and Theorem \ref{thm:levyExponent} are the main results of our paper.
Let us recall definitions of $\tau_{U}$ and $\tau_{D}$ stated before
Proposition \ref{prop:skorohodEmbending} and that $g^{c}$ is the
unique minimizer in the problem \eqref{eq:variational}. These quantities
will be used as a functions of Brownian paths. We define the processes
$\cbr{X_{t}^{c}}_{t\in\R}$ and $\cbr{Y_{t}^{c}}_{t\in\R}$ by putting
$X_{t}^{c}=Y_{t}^{c}=0$ for $t\leq0$ and 
\begin{align*}
X_{t}^{c} & \eqdef\sum_{k=-\infty}^{+\infty}L_{t}^{2kc}-\sum_{k=-\infty}^{+\infty}L_{t}^{(2k+1)c},\\
Y_{t}^{c} & \eqdef\sum_{k=-\infty}^{+\infty}L_{t}^{kc},
\end{align*}
where $L^{a}$ is the local time of a Brownian motion at the level
$a\in\R$. 
\begin{thm}
\label{thm:localTimeRepresentation}Let $B$ be a Brownian motion
independent of $X^{c},Y^{c}$. Then 
\begin{equation}
B_{t}^{c}=^{d}\begin{cases}
\rbr{B_{\tau_{U}}-c/2}+X_{t-\tau_{U}}^{c}, & \text{ if }\tau_{U}<\tau_{D},\\
\rbr{B_{\tau_{D}}+c/2}-X_{t-\tau_{D}}^{c}, & \text{ if }\tau_{U}>\tau_{D}.
\end{cases}\label{eq:rep}
\end{equation}
and 
\begin{equation}
\TTV B{[0;t]}c=\text{TV}(B^{c})[0,t]=^{d}Y_{t-(\tau_{U}\wedge\tau_{D})}^{c}.\label{eq:representtation}
\end{equation}
In both cases above $=^{d}$, the equality in distribution, is understood
on the process level. \end{thm}
\begin{rem}
It is easy to prove that $\tau_{U}\wedge\tau_{D}$ has exponentially
decaying tails. Further, we note that $\tau_{U}\wedge\tau_{D}$ is
the same as $\theta(c)$ in \cite{Imhof:1985fk}. Thus its density
is given by \cite[(2.3)]{Imhof:1985fk}. The Laplace transform is
also provided therein. 
\end{rem}
The process $Y^{c}$ is non-decreasing and can be given a straightforward
description in terms of Lévy processes. We define its generalized
inverse $\cbr{S_{t}^{c}}_{t\geq0}$ by 
\[
S_{t}^{c}\eqdef\inf\cbr{s\geq0:Y_{s}^{c}=t}.
\]

\begin{thm}
\label{thm:levyExponent}The process $S^{c}$ is a Lévy process with
the exponent \textup{$\Phi^{c}(q)\eqdef(-\log\ev{e^{-qS_{t}^{c}}})/t$,
given by} 
\begin{equation}
\Phi^{c}(q)=\sqrt{2q}\tanh\rbr{c\sqrt{q/2}},\quad q\geq0.\label{eq:exponent}
\end{equation}
\end{thm}
\begin{rem}
\label{rem:comparingToOldResults}Potentially this result can be used
to reprove \cite[Theorem 5]{Lochowski:2013yq} in case of $X$ being
a Brownian motion. Indeed, the process $\cbr{c^{-1}S_{t}^{c}}_{t\geq0}$
is inverse of the process $c\:\TTV Xtc$. By Theorem \ref{thm:levyExponent}
it is a Lévy process with the exponent $q\mapsto\Phi^{c}(c^{-1}q)$. 
\end{rem}

\begin{rem}
Without loss of generality we may choose $c=1/2$ then the \eqref{eq:representtation}
can be regarded as the local time of a Brownian motion on circle at
point $1/2$. Processes of this kind were studied e.g. in \cite{Pitman:1996uq,Bolthausen:1979fk}.
\end{rem}

\begin{rem}
The theorem gives also a representation for the two-sided Skorokhod
problem for Brownian motion, like the one studied in \cite{Kruk:2007fk}.\end{rem}
\begin{proof}
(of Theorem \ref{thm:localTimeRepresentation}). By Proposition \ref{prop:skorohodEmbending}
for $t\leq\tau_{U}\wedge\tau_{D}$ the process $B^{c}$ is constant.
For the rest of the proof we assume, without loss of generality, that
$\tau_{U}<\tau_{D}$ (the other case follow by taking $-B$). We have
$x=B_{\tau_{U}}-c/2$ and we want to identify the law joint low of
$f_{1}$ and $f_{2}$ (the path-wise analogues of $f_{1}$ and $ $$f_{2}$
in Proposition \ref{prop:skorohodEmbending}). We have $f_{1}(t)=f_{2}(t)=0$
for $t\in[0,\tau_{U}]$. At the terminal point of this interval we
have 
\[
B_{\tau_{U}}-(x+f_{1}(\tau_{U})+f_{2}(\tau_{U}))=c/2.
\]
By the unicity of solutions of the Skorokhod problem the task boils
down to ensuring that $[0,+\infty)\ni t\mapsto f_{1}(\tau_{U}+t)=:\tilde{f}_{1}(t)$
and $[0,+\infty)\ni t\mapsto f_{2}(\tau_{U}+t)=:\tilde{f}_{2}(t)$
fulfills conditions of Lemma \ref{lem:skorohod} with for $x=-c/2$
and the function $t\mapsto B_{\tau_{u}+t}-B_{\tau_{U}}$. By the strong
Markov property this process is a Brownian motion independent of its
evolution up to time $\tau_{U}$. 

We are going to find the joint law of $\tilde{f}_{1}$ and $\tilde{f}_{2}$.
To this end we consider the function $F_{c}$ given by
\[
F_{c}(x):=\sum_{k=-\infty}^{+\infty}|x-2kc|\mathbf{1}_{x\in[(2k-1)c,(2k+1)c)}-c/2.
\]
One easily calculates its left derivative
\[
\partial^{-}F_{c}(x)=\sum_{k=-\infty}^{+\infty}\mathbf{1}_{x\in(2kc,(2k+1)c]}-\sum_{k=-\infty}^{+\infty}\mathbf{1}_{x\in((2k-1)c,2kc]}.
\]
and the weak second derivative 
\[
\partial_{2}F_{c}=\sum_{k=-\infty}^{+\infty}\delta_{2kc}-\delta_{(2k+1)c}.
\]
Let $B$ be a Brownian motion. Let us consider $F_{c}(B)$. By the
It\={o}-Tanaka formula \cite[Theorem VI.1.5]{Revuz:1991kx} (the formula
as stated applies only to convex functions but it is easy to represent
$F_{c}=F_{c}^{1}-F_{c}^{2}$ where both $F_{c}^{1}$ and $F_{c}^{2}$
are piecewise linear and convex) we have
\begin{equation}
F_{c}(B_{t})=\int_{0}^{t}\partial^{-}F_{c}(B_{s})\dd B_{s}+\sum_{k=-\infty}^{+\infty}L_{t}^{2kc}-\sum_{k=-\infty}^{+\infty}L_{t}^{(2k+1)c}.\label{eq:basic}
\end{equation}
Let us notice that by the Lévy theorem $\beta_{t}:=\int_{0}^{t}\partial^{-}F_{c}(B_{s})\dd B_{s}$
defines another Brownian motion. A path-wise application of Lemma
\ref{lem:skorohod} reveals that $x=c/2$, $h(t)=F_{c}(B_{t})$, $f(t)=\beta_{t}$,
$\tilde{f}_{1}(t)=\sum_{k=-\infty}^{+\infty}L_{t}^{2kc},$ $\tilde{f_{2}}(t)=\sum_{k=-\infty}^{+\infty}L_{t}^{(2k+1)c}$
is the solution of the Skorokhod problem. Putting the pieces together
we obtain the representation \eqref{eq:rep}.
\end{proof}

\begin{proof}
(of Theorem \ref{thm:levyExponent}) First, we are going to prove
that $S^{c}$ is a Lévy process. We will check conditions of \cite[Definition 1.6]{Sato:1999fk}.
Points 2 and 5 are straightforward. Let us prove 1 for the case $0\leq t_{0}<t_{1}$
(the case of general $n$ follows by induction). $S_{t_{1}}^{c}$
is clearly a stopping time with respect to the filtration of the underlying
Brownian motion. We have $S_{t_{2}}^{c}-S_{t_{1}}^{c}=\inf\cbr{s\geq0:Y_{s+S_{t_{1}}^{c}}^{c}-Y_{S_{t_{1}}^{c}}^{c}\geq t_{2}-t_{1}}=\inf\cbr{s\geq0:\theta_{S_{t_{1}}^{c}}\circ Y_{s}^{c}\geq t_{2}-t_{1}}$,
where $\theta$ is the shift operator. By the strong Markov property
we conclude that $S_{t_{2}}^{c}-S_{t_{1}}^{c}$ is independent of
$S_{t_{1}}^{c}$ concluding the proof of 1. A very similar argument
covers 3. Finally 4, can be proved by the fact that $\pr{L_{t}^{0}>0}=1$
for any $t>0$. 

Further, the proof will use many notions of the theory of Brownian
motion. As they are standard, instead of introducing them formally
(which would be very lengthy), we refer the reader to \cite{Revuz:1991kx}.

Process $S^{c}$ is clearly non-decreasing (i.e. it is a subordinator).
Our aim is to calculate its Lévy exponent $\Phi^{c}$. To this end
let us denote by $n$ the It\={o} excursion measure (for details we
refer to \cite[Chapter XII]{Revuz:1991kx}). Let $\zeta$ denote the
length of the excursion and $ $$\rho=\inf\cbr{s\geq0:|w_{s}|=c}$
(by $w$ we denote the excursion itself); by convention we put $\rho=+\infty$
if the defining set is empty. By the It\={o} decomposition, \cite[Theorem XII.2.4]{Revuz:1991kx}
for any $q\geq0$ we have 

\begin{equation}
\Phi^{c}(q)=\int_{(0,\infty)}(1-e^{-qy})\rbr{n(\rho\in\dd y;\rho<+\infty)+n(\zeta\in\dd y;\rho=+\infty)}=:I_{1}+I_{2}.\label{eq:tmp1}
\end{equation}
Let us now denote the stopping time $\tau_{a}\eqdef\inf\cbr{s\geq0:B_{s}=a}$
for $a\in\R$ and $\mathbb{P}_{x}$ the measure under which the Brownian
motion $B$ starts from $x$. \cite[Proposition 2]{Doney:2005fk}
suggests that for some constant $k>0$ we have 
\begin{equation}
I_{1}=\lim_{x\searrow0}I_{1}(x),\quad I_{1}(x)\eqdef kx^{-1}\int_{(0,\infty)}(1-e^{-qy})\mathbb{P}_{x}(\tau_{c}\in\dd y;\tau_{c}<\tau_{0}),\label{eq:exp1}
\end{equation}
and 
\begin{equation}
I_{2}=\lim_{x\searrow0}I_{2}(x),\quad I_{2}(x)\eqdef kx^{-1}\int_{(0,\infty)}(1-e^{-qy})\mathbb{P}_{x}(\tau_{0}\in\dd y;\tau_{c}>\tau_{0}).\label{eq:exp2}
\end{equation}
Proving these relations is surprisingly lengthy so we postpone it
until later. Now we are going to show how they imply our result. By
\cite[(8.8)]{Kyprianou:2006fk} we have
\[
I_{1}=k\lim_{x\searrow0}x^{-1}\mathbb{E}_{x}(1-e^{-q\tau_{c}};\tau_{c}<\tau_{0})=k\lim_{x\searrow0}x^{-1}\rbr{\frac{W^{(0)}(x)}{W^{(0)}(c)}-\frac{W^{(q)}(x)}{W^{(q)}(c)}},
\]
where $W^{(q)}(x)=\sqrt{2/q}\sinh(\sqrt{2q}x)$ and $W^{(0)}=2x$
(as indicated in \cite[Excercise 8.2 p.233]{Kyprianou:2006fk}). Simple
calculations reveal that 
\[
I_{1}=k\rbr{\frac{1}{c}-\frac{\sqrt{2q}}{\sinh(\sqrt{2q}c)}}.
\]
Similarly, by \cite[(8.9)]{Kyprianou:2006fk}, we have
\[
I_{2}=k\lim_{x\searrow0}x^{-1}\mathbb{E}_{x}(1-e^{-q\tau_{0}};\tau_{0}<\tau_{c})=k\lim_{x\searrow0}x^{-1}\rbr{Z^{(0)}(x)-Z^{(0)}(c)\frac{W^{(0)}(x)}{W^{(0)}(c)}-Z^{(q)}(x)+Z^{(q)}(c)\frac{W^{(q)}(x)}{W^{(q)}(c)}},
\]
where $Z^{(q)}(x)=\cosh(\sqrt{2q}x)$ and $Z^{(0)}(x)=1$. Performing
some standard calculations we get
\[
I_{2}=k\rbr{\sqrt{2q}\coth(\sqrt{2q}c)-\frac{1}{c}}.
\]
We thus have $I_{1}+I_{2}=k\sqrt{2q}\tanh(c\sqrt{q/2})$. The constant
$k$ have yet to be determined. This, in principle could be done using
\cite{Doney:2005fk} but we shall do this by comparing to results
of \cite[Theorem 1]{Lochowski:2013yq} in the case of a Brownian motion.
As indicated in Remark \ref{rem:comparingToOldResults} we are to
study a Lévy process with the exponent $q\mapsto\Phi^{c}(c^{-1}q)$.
One easily checks that $\lim_{c\searrow0}\Phi^{c}(c^{-1}q)=kq$ which
describes a drift process with speed $2k$. By $ $\cite[Theorem 1]{Lochowski:2013yq}
we conclude that $k=1$. To avoid unnecessary notation we will omit
writing $k$ in the further part of the proof. 

The last step of this proof is justifying \eqref{eq:exp1} and \eqref{eq:exp2}.
In both the cases we will introduce two additional quantities $I_{i}^{m}$
and $I_{i}^{m}(x)$, where $i\in\cbr{1,2}$. The parameter $m$ controls
some discretisation which is required to apply \cite[Proposition 2]{Doney:2005fk}.
We will show that for $i\in\cbr{1,2}$ we have 
\begin{equation}
\lim_{m\to+\infty}I_{i}^{m}=I_{i},\label{eq:convergence1}
\end{equation}
\begin{equation}
\lim_{m\to+\infty}\sup_{x\in(0,1/2)}|I_{i}^{m}(x)-I_{i}(x)|=0,\label{eq:convergence2}
\end{equation}
\begin{equation}
\forall_{m}\lim_{x\searrow0}I_{i}^{m}(x)=I_{i}^{m}.\label{eq:convergence3}
\end{equation}
This will be enough to show \eqref{eq:exp1} and \eqref{eq:exp2}.
Indeed, we know already that $\lim_{x\searrow0}I_{i}(x)$ exists,
it is enough to find a sequence $\cbr{x_{n}}_{n\geq0}$ such that
$x_{n}\to0$ and $I_{i}(x_{n})\to I_{i}$. To this end we fix $\epsilon>0$
and choose $m$ such that $|I_{i}^{m}-I_{i}|\leq\epsilon$ and $\sup_{x\in(0,1/2)}|I_{i}^{m}(x)-I_{i}(x)|\leq\epsilon$,
finally we can find $x\in(0,1/2)$ such that $|I_{i}^{m}(x)-I_{i}^{m}|\leq\epsilon$.
For this $x$ we have $|I_{i}-I_{i}(x)|\leq3\epsilon$. Further reasoning
is standard. 

Let us start with \eqref{eq:exp2}. We define $\cbr{a_{k}^{m}}_{k\in\cbr{-\infty,\ldots,+\infty}}$
where $m\in\mathbb{N}$ by 
\begin{equation}
a_{k}^{m}\eqdef\begin{cases}
\rbr{1-\frac{k}{m}}^{-2} & \text{ for }k\leq0,\\
1+\frac{k}{m} & \text{ for }k>0.
\end{cases}\label{eq:akm}
\end{equation}
One easily verifies that for any $m$ we have $a_{k}^{m}\leq a_{k+1}^{m}$.
By the fact that function $x\mapsto1-e^{-qx}$ is increasing we have
\[
J^{m}\eqdef\sum_{k=-\infty}^{+\infty}\rbr{1-e^{-qa_{k-1}^{m}}}n(\zeta\in[a_{k-1}^{m},a_{k}^{m});\rho=+\infty)\leq I_{2}\leq\sum_{k=-\infty}^{+\infty}\rbr{1-e^{-qa_{k}^{m}}}n(\zeta\in[a_{k-1}^{m},a_{k}^{m});\rho=+\infty)=:K^{m}.
\]
We are going to prove that $\lim_{m\to+\infty}K^{m}-J^{m}=0$ and
consequently 
\begin{equation}
\lim_{m\to+\infty}K^{m}=I_{2}.\label{eq:convergence10}
\end{equation}
We have%
\footnote{From now on notation $x\cleq y$ denotes a situation when there exists
a constant $C>0$ such that $x\leq Cy$ and $C$ is irrelevant for
calculations.%
}
\begin{align}
0\leq K^{m}-J^{m} & \leq\sum_{k=-\infty}^{+\infty}\rbr{e^{-qa_{k-1}^{m}}-e^{-qa_{k}^{m}}}n(\zeta\in[a_{k-1}^{m},a_{k}^{m}))\label{eq:krowaZKrakowa}\\
 & \cleq\sum_{k=-\infty}^{0}(a_{k}^{m}-a_{k-1}^{m})n(\zeta\in[a_{k-1}^{m},a_{k}^{m}))+\frac{1}{m}\sum_{k=1}^{+\infty}e^{-qa_{k-1}^{m}}n(\zeta\in[a_{k-1}^{m},a_{k}^{m})),\nonumber 
\end{align}
where we used $\rbr{e^{-qa_{k-1}^{m}}-e^{-qa_{k}^{m}}}\cleq e^{-qa_{k-1}^{m}}(a_{k}^{m}-a_{k-1}^{m})$.
By the special choice of sequence $a^{m}$ and \cite[Proposition XII.2.8]{Revuz:1991kx}
for some $C>0$ we have $n(\zeta\in(a_{k-1}^{m},a_{k}^{m}])\leq C/m$
which holds for $k\leq0$. The first sum can be upper-bounded by 
\[
\frac{1}{m}\sum_{k=-\infty}^{0}(a_{k}^{m}-a_{k-1}^{m})=\frac{1}{m}a_{0}^{m}\to0.
\]
The second term of \eqref{eq:krowaZKrakowa} is bounded from above
by the following integral 
\[
\frac{1}{m}\int_{1}^{+\infty}e^{-qy}n(\zeta\in\dd y)\to0.
\]
We have thus proven \eqref{eq:convergence10}. Now we define the aforementioned
$I_{2}^{m}$, namely we put 
\[
I_{2}^{m}\eqdef\sum_{k=-\infty}^{+\infty}\rbr{1-e^{-qa_{k}^{m}}}n(\zeta\in[a_{k-1}^{m},a_{k}^{m});a_{k-1}^{m}\leq\rho).
\]
Obviously, we have $K^{m}\leq I_{2}^{m}$ and further we can estimate
\begin{align*}
0\le I_{2}^{m}-K^{m} & \leq\sum_{k=-\infty}^{+\infty}\rbr{1-e^{-qa_{k}^{m}}}n(\zeta\in[a_{k-1}^{m},a_{k}^{m});\rho\in[a_{k-1}^{m},a_{k}^{m}])\\
 & \leq\int\rbr{1-e^{-q\zeta(w)}}\sum_{k=-\infty}^{+\infty}1_{\zeta(w)\in[a_{k-1}^{m},a_{k}^{m})}1_{\rho(w)\in[a_{k-1}^{m},a_{k}^{m}]}n(\dd w).
\end{align*}
In the last expression we integrate over the space of excursions (we
refer the reader to \cite[Section XII]{Revuz:1991kx} for details).
By \cite[Proposition XII.2.8]{Revuz:1991kx} one checks that $\int\rbr{1-e^{-q\zeta(w)}}n(\dd w)<+\infty$.
The expression 
\[
\sum_{k=-\infty}^{+\infty}1_{\zeta(w)\in(a_{k-1}^{m},a_{k}^{m}]}1_{\rho(w)\in[a_{k-1}^{m},a_{k}^{m}]}
\]
 is bounded by $1$ and converges point-wise to $0$ (for any $w$).
Thus by Lebesgue's dominated convergence we get that $\lim_{m\to+\infty}I_{2}^{m}-K^{m}.$
This and \eqref{eq:convergence10} yields \eqref{eq:convergence1}.

Now we define the aforementioned $I_{2}^{m}(x)$ by 
\begin{equation}
I_{2}^{m}(x)\eqdef\sum_{k=-\infty}^{+\infty}\rbr{1-e^{-qa_{k}^{m}}}\frac{1}{x}\mathbb{P}_{x}(\tau_{0}\in[a_{k-1}^{m},a_{k}^{m});a_{k-1}^{m}\leq\tau_{c}).\label{eq:i2mx}
\end{equation}
We note that by \cite[Proposition 2]{Doney:2005fk} every summand
converge to the corresponding summand of $I_{2}^{m}$. Let us consider
terms with $k\geq1$. We have
\begin{align*}
\frac{1}{x}\mathbb{P}_{x}(\tau_{0}\in[a_{k-1}^{m},a_{k}^{m}];a_{k-1}^{m}\leq\tau_{c}) & \leq\frac{\mathbb{P}_{x}(\tau_{0}\geq1)}{x}\mathbb{P}_{x}(a_{k-1}^{m}\leq\tau_{c};a_{k-1}^{m}\leq\tau_{0}|\tau_{0}\geq1)\\
 & \cleq\exp\rbr{-Ca_{k-1}^{m}},
\end{align*}
for some constant $C>0$. Indeed, $\mathbb{P}_{x}(\tau_{0}\geq1)=\mathbb{P}_{0}(\inf_{t\in[0,1]}W_{t}\geq-x)\approx(2/\pi)^{1/2}x$
(see \cite[Section III.3.7]{Revuz:1991kx}). Secondly, by the strong
Markov property $\mathbb{P}_{x}(a_{k-1}^{m}\leq\tau_{c};a_{k-1}^{m}\leq\tau_{0}|\tau_{0}\geq1)\leq\mathbb{P}_{c/2}(\forall_{t\in[0,a_{k-1}^{m}-1]}W_{t}\in(0,c))\cleq\exp\rbr{-C(a_{k-1}^{m}-1)}$.
The last term is clearly summable, thus Lebesgue's dominated theorem
implies that 
\[
\lim_{x\searrow0}\sum_{k=1}^{+\infty}\frac{1}{x}\mathbb{P}_{x}(\tau_{0}\in[a_{k-1}^{m},a_{k}^{m}];a_{k-1}^{m}\leq\tau_{c})=\sum_{k=1}^{+\infty}\rbr{1-e^{-qa_{k}^{m}}}n(\zeta\in[a_{k-1}^{m},a_{k}^{m}];a_{k-1}^{m}\leq\rho).
\]
Next, we treat the case of $k\leq0$. Using \cite[Remark 1 after Proposition III.3.8]{Revuz:1991kx}
we estimate
\begin{align}
\rbr{1-e^{-qa_{k}^{m}}}\frac{1}{x}\mathbb{P}_{x}(\tau_{0}\in[a_{k-1}^{m},a_{k}^{m}];a_{k-1}^{m}\leq\tau_{c}) & \leq\rbr{1-e^{-qa_{k}^{m}}}\frac{1}{x}\int_{a_{k-1}^{m}}^{a_{k}^{m}}\frac{x}{(2\pi y^{3})^{1/2}}\exp\rbr{-\frac{x^{2}}{2y}}\dd y\label{eq:tralala}\\
 & \leq\rbr{1-e^{-qa_{k}^{m}}}\frac{a_{k}^{m}-a_{k-1}^{m}}{(a_{k-1}^{m})^{3/2}}\leq\frac{a_{k}^{m}-a_{k-1}^{m}}{(a_{k-1}^{m})^{1/2}}.\nonumber 
\end{align}
One checks that $a_{k}^{m}-a_{k-1}^{m}\cleq(1-k/m)^{-3}$, this implies
that the last expression is a sequence summable in $k$. Analogously
as before we have convergence for $\sum_{k=-\infty}^{0}\ldots$ .
Put together they give \eqref{eq:convergence3}.

We notice that $I_{2}^{m}(x)$ defined in \eqref{eq:i2mx} can be
expressed as 
\begin{equation}
I_{2}^{m}(x)=\int_{(0,+\infty)}\rbr{1-e^{-qf^{m}(y)}}\frac{1}{x}\mathbb{P}_{x}(\tau_{0}\in\dd y;g^{m}(y)\leq\tau_{c}),\label{eq:notation}
\end{equation}
where $g^{m}(y)=a_{k-1}^{m}$ and $f^{m}(y)=a_{k}^{m}$ whenever $y\in[a_{k-1}^{m},a_{k}^{m})$.
We recall $I_{2}(x)$ defined in \eqref{eq:exp2} and consider 
\begin{align*}
|I_{2}^{m}(x)-I_{2}(x)| & \leq\int_{(0,+\infty)}\rbr{e^{-qy}-e^{-qf^{m}(y)}}\frac{1}{x}\mathbb{P}_{x}(\tau_{0}\in\dd y;g^{m}(y)\leq\tau_{c})\\
 & +\int_{(0,+\infty)}\rbr{1-e^{-qy}}\frac{1}{x}\mathbb{P}_{x}(\tau_{0}\in\dd y;\tau_{c}\in[g^{m}(y),y])=:J^{m}(x)+K^{m}(x).
\end{align*}
We estimate 
\[
J^{m}(x)\cleq\int_{(0,1)}\frac{e^{-qy}-e^{-qf^{m}(y)}}{y}\frac{y}{x}\mathbb{P}_{x}(\tau_{0}\in\dd y)+\int_{[1,+\infty)}\rbr{e^{-qy}-e^{-qf^{m}(y)}}\mathbb{P}_{x}(\tau_{c}\geq g^{m}(y)|\tau_{0}\geq1)\dd y.
\]
We denote $l^{m}(y)\eqdef\frac{e^{-qy}-e^{-qf^{m}(y)}}{y}$. One checks
that $\forall_{y\in(0,1]}$ we have $l^{m}(y)\leq2$ and $\lim_{m\to+\infty}l^{m}(y)=0$.
For $y\geq1$ we have $e^{-qy}-e^{-qf^{m}(y)}\cleq1/m$ and $\mathbb{P}_{x}(\tau_{c}\geq g^{m}(y)|\tau_{0}\geq1)\cleq e^{-Cy}$.
We recall also \cite[Remark 1 after Proposition III.3.8]{Revuz:1991kx}
to get 
\[
J^{m}(x)\cleq\int_{(0,1)}l^{m}(y)\frac{1}{(2\pi y)^{1/2}}\exp\rbr{-\frac{x^{2}}{2y}}\dd y+\frac{1}{m}\int_{[1,+\infty)}e^{-Cy}\dd y\cleq\int_{(0,1)}\frac{l^{m}(y)}{(2\pi y)^{1/2}}\dd y+\frac{1}{m}.
\]
We can see that the last expression does not depend on $x$ and converges
to $0$ when $m\to+\infty$. We proceed to $K^{m}$. We recall that
$\mathbb{P}_{x}(\tau_{c}\leq\tau_{0})\cleq x$ and thus
\[
K^{m}(x)\cleq\mathbb{P}_{x}(\tau_{c}\in(g^{m}(\tau_{0}),\tau_{0})|\tau_{c}\leq\tau_{0}).
\]
We observe that the strong Markov property yields $\mathbb{P}_{x}(\tau_{c}\in(g^{m}(\tau_{0}),\tau_{0})|\tau_{c}\leq\tau_{0})\leq\mathbb{P}_{c}(\tau_{0}\leq1/m)\to0$.
This is enough to conclude that $\sup_{x\in(0,1)}K^{m}(x)\to0$ when
$m\to+\infty$. And thus also \eqref{eq:convergence2}. Having checked
that \eqref{eq:convergence1}, \eqref{eq:convergence2} and \eqref{eq:convergence3}
hold for $i=2$ we conclude \eqref{eq:exp2}. 

Now we turn to \eqref{eq:exp1}. We use $\cbr{a_{k}^{m}}_{k\in\cbr{-\infty,\ldots,+\infty}}$
given by \eqref{eq:akm} and define 
\[
I_{1}^{m}\eqdef\sum_{l}\sum_{k>l}\rbr{1-e^{-qa_{l}^{m}}}n(\rho\in[a_{l}^{m},a_{l+1}^{m});\zeta\in[a_{k}^{m},a_{k+1}^{m})).
\]
Obviously we have $I_{1}^{m}\leq I_{1}$ and further 
\begin{align*}
I_{1}-I_{1}^{m} & \leq\sum_{l}\rbr{e^{-qa_{l}^{m}}-e^{-qa_{l+1}^{m}}}n(\rho\in[a_{l}^{m},a_{l+1}^{m});\rho<+\infty)\\
 & +\sum_{k}\rbr{1-e^{-qa_{k}^{m}}}n(\rho\in[a_{k}^{m},a_{k+1}^{m});\zeta\in[a_{k}^{m},a_{k+1}^{m}))=:K^{m}+J^{m}.
\end{align*}
We consider $K^{m}$. We have
\[
K^{m}=\int l^{m}(y)(1-e^{-qy})n(\rho\in\dd y;\rho<+\infty),
\]
where $l^{m}(y)\eqdef\rbr{e^{-qf^{m}(y)}-e^{-qg^{m}(y)}}/(1-e^{-qy})$.
We have 
\begin{equation}
\int_{(0,+\infty)}(1-e^{-qy})n(\rho\in\dd y;\rho<+\infty)=\int_{(0,+\infty)}qe^{-qy}n(\rho\geq y;\rho<+\infty)\dd y<+\infty,\label{eq:finiteAAA}
\end{equation}
where the last estimate follows by $n(\rho\geq y;\rho<+\infty)\leq n(\zeta>y)\cleq y^{-1/2}$,
which is asserted in \cite[Proposition XII.2.8]{Revuz:1991kx}. It
is also easy to verify that $\sup_{y>0}l^{m}(y)\leq C$ for some $C>0$
and that $l^{m}\to0$ point-wise. Dominated Lebesgue's theorem implies
$K^{m}\to0$. To deal with $J^{m}$ one checks that $1-e^{-qg^{m}(y)}\cleq1-e^{-qy}$
\begin{align*}
J^{m}\cleq & \int_{(0,+\infty)}\rbr{1-e^{-qy}}n(\rho\in\dd y;\zeta\in[g^{m}(y),f^{m}(y))\dd y\\
\leq & \int_{(0,\epsilon)}\rbr{1-e^{-qy}}n(\rho\in\dd y)\dd y+\int_{(\epsilon,+\infty)}n(\rho\in\dd y;\zeta\in[g^{m}(y),f^{m}(y))\dd y.
\end{align*}
The measure $n$ is finite on the set $\rho\geq\epsilon$ hence the
second term converges by the Lebesgue dominated theorem (as the conditions
converge to $0$). Further, by \eqref{eq:finiteAAA}, the first integral
is equal to some $C(\epsilon)>0$, which $\lim_{\epsilon\to0}C(\epsilon)=0$.
The above facts are enough to conclude that $J^{m}\to0$. In this
way we have established \eqref{eq:convergence1}. 

Now we define the approximation sequence $I_{1}^{m}(x)$ by
\begin{equation}
I_{1}^{m}(x)\eqdef\sum_{l}\sum_{k>l}\rbr{1-e^{-qa_{l}^{m}}}\frac{1}{x}\mathbb{P}_{x}\rbr{\tau_{c}\in[a_{l}^{m},a_{l+1}^{m});\tau_{0}\in[a_{k}^{m},a_{k+1}^{m})}.\label{eq:I1mx}
\end{equation}
We are going to show \eqref{eq:convergence3}. To this end we define
\[
I_{1}^{m}(x;h)=\sum_{l=-h}^{h}\sum_{k=l+1}^{h}\rbr{1-e^{-qa_{l}^{m}}}\frac{1}{x}\mathbb{P}_{x}\rbr{\tau_{c}\in[a_{l}^{m},a_{l+1}^{m});\tau_{0}\in[a_{k}^{m},a_{k+1}^{m})},\quad h\in\cbr{2,3,\ldots}.
\]
Each of $I_{1}^{m}(x;h)$ contains only finite number of terms so
by \cite[Proposition 2]{Doney:2005fk} we have 
\[
\lim_{x\searrow0}I_{1}^{m}(x;h)=\sum_{l=-h}^{h}\sum_{k=l+1}^{h}\rbr{1-e^{-qa_{l}^{m}}}n\rbr{\tau_{c}\in[a_{l}^{m},a_{l+1}^{m});\tau_{0}\in[a_{k}^{m},a_{k+1}^{m})}.
\]
To conclude \eqref{eq:convergence3} we are going to show that for
any fixed $m$ we have 
\begin{equation}
\lim_{h\to+\infty}\sup_{x}(I_{1}^{m}(x)-I_{1}^{m}(x;h))=0.\label{eq:uciete}
\end{equation}
Firstly, we notice that 
\[
\sum_{k=h}^{+\infty}\sum_{l\leq k}\rbr{1-e^{-qa_{l}^{m}}}\frac{1}{x}\mathbb{P}_{x}\rbr{\tau_{c}\in[a_{l}^{m},a_{l+1}^{m});\tau_{0}\in[a_{k}^{m},a_{k+1}^{m})}\leq\sum_{k=h}^{+\infty}\sum_{l\leq k}\frac{1}{x}\mathbb{P}_{x}\rbr{\tau_{0}\in[a_{k}^{m},a_{k+1}^{m})}\leq\frac{\mathbb{P}_{x}\rbr{\tau_{0}\geq a_{h}^{m}}}{x}.
\]
The last term converges to $0$ with $h\to+\infty$ uniformly in $x$.
Next, we treat
\begin{align*}
\sum_{l\leq-h}\sum_{k=l+1}^{h}\rbr{1-e^{-qa_{l}^{m}}}\frac{1}{x}\mathbb{P}_{x}\rbr{\tau_{c}\in[a_{l}^{m},a_{l+1}^{m});\tau_{0}\in[a_{k}^{m},a_{k+1}^{m})} & \leq\sum_{l\leq-h}\frac{1}{x}\mathbb{P}_{x}\rbr{\tau_{c}\in[a_{l}^{m},a_{l+1}^{m});\tau_{0}\geq\tau_{c}}\\
 & \leq\frac{1}{x}\mathbb{P}_{x}\rbr{\tau_{c}\leq a_{-h}^{m};\tau_{0}\geq\tau_{c}}.
\end{align*}
We are now to analyze the last expression. By the strong Markov property
we have
\begin{align*}
\frac{1}{x}\mathbb{P}_{x}\rbr{\tau_{c}\leq a_{-h}^{m};\tau_{0}\geq\tau_{c}} & \leq\frac{1}{x}\mathbb{P}_{x}\rbr{\tau_{c}\leq a_{-h}^{m};\tau_{0}\geq\tau_{c/2}}=\frac{\mathbb{P}_{x}(\tau_{0}\geq\tau_{c/2})}{x}\mathbb{P}_{x}\rbr{\tau_{c}\leq a_{-h}^{m}|\tau_{0}\geq\tau_{c/2}}\\
 & =\frac{\mathbb{P}_{x}(\tau_{0}\geq\tau_{c/2})}{x}\mathbb{P}_{c/2}(\tau_{c}\leq a_{-h}^{m}-\tau_{c/2})\cleq\mathbb{P}_{c/2}(\tau_{c}\leq a_{-h}^{m})\conv0.
\end{align*}
We notice that the last expression does not involve $x$, thus we
again we have obtained convergence uniform in $x$ and finish the
proof of \eqref{eq:uciete} and consequently \eqref{eq:convergence3}. 

Our final task is \eqref{eq:convergence2}. We recall \eqref{eq:exp1}
and we decompose 
\begin{align*}
I_{1}(x)= & \sum_{l}\sum_{k>l}\int_{[a_{l}^{m},a_{l+1}^{m})}\rbr{1-e^{-qy}}\frac{1}{x}\mathbb{P}_{x}\rbr{\tau_{c}\in\dd y;\tau_{0}\in[a_{k}^{m},a_{k+1}^{m})}\\
+ & \sum_{k}\int_{[a_{k}^{m},a_{k+1}^{m})}\rbr{1-e^{-qy}}\frac{1}{x}\mathbb{P}_{x}\rbr{\tau_{c}\in\dd y;\tau_{c}\leq\tau_{0};\tau_{0}\in[a_{k}^{m},a_{k+1}^{m})}.
\end{align*}
We recall \eqref{eq:I1mx} and consider
\begin{align*}
|I_{1}^{m}(x)-I_{1}(x)| & \cleq\sum_{l}\sum_{k>l}\int_{[a_{l}^{m},a_{l+1}^{m})}\rbr{e^{-qa_{l}^{m}}-e^{-qa_{l+1}^{m}}}\frac{1}{x}\mathbb{P}_{x}\rbr{\tau_{c}\in\dd y;\tau_{0}\in[a_{k}^{m},a_{k+1}^{m})}\\
 & +\sum_{k}\rbr{1-e^{-qa_{k+1}^{m}}}\frac{1}{x}\mathbb{P}_{x}\rbr{\tau_{c}\in[a_{k}^{m},a_{k+1}^{m});\tau_{c}\leq\tau_{0};\tau_{0}\in[a_{k}^{m},a_{k+1}^{m})}=:J^{m}(x)+K^{m}(x).
\end{align*}
We deal with the first term as follows 
\begin{align*}
J^{m}(x) & \cleq\sum_{l}\sum_{k>l}\int_{[a_{l}^{m},a_{l+1}^{m})}(a_{l+1}^{m}-a_{l}^{m})\frac{1}{x}\mathbb{P}_{x}\rbr{\tau_{c}\in\dd y;\tau_{0}\in(a_{k}^{m},a_{k+1}^{m}];\tau_{c/2}\leq\tau_{0}}\\
 & \cleq\sum_{l}\sum_{k>l}\int_{[a_{l}^{m},a_{l+1}^{m})}(a_{l+1}^{m}-a_{l}^{m})\mathbb{P}_{x}\rbr{\tau_{c}-\tau_{c/2}\in\dd y;\tau_{0}+\tau_{c/2}\in(a_{k}^{m},a_{k+1}^{m}]|\tau_{c/2}\leq\tau_{0}}\\
 & \cleq\max_{l}(a_{l+1}^{m}-a_{l}^{m}).\\
\end{align*}
Obviously the last expression independent of $x$ and convergences
to $0$. Similarly, one can estimate 
\[
K^{m}(x)\cleq\sum_{k}\mathbb{P}_{c/2}\rbr{\tau_{c}\in(a_{k}^{m},a_{k+1}^{m}];\tau_{0}\in(a_{k}^{m},a_{k+1}^{m}]}.
\]
It is easy to see that this converge to $0$. In this way we have
shown \eqref{eq:convergence2}. This completes the proof. 
\end{proof}

\subsection*{Acknowledgments }

The author wish to thank prof. Andreas Kyprianou for useful discussions
and in particular for the ideas which lead to the proof of Theorem
\ref{thm:levyExponent}. Further, the author is grateful to dr. R.
\L{}ochowski for useful comments and suggestions. 

The research was partially supported by MNiSW grant N N201 397537.

\bibliographystyle{abbrv}
\bibliography{/Users/piotrmilos/Dropbox/Sync/Library/branching}

\end{document}